\documentclass[12pt]{amsart}

\newcommand{\bfA}{\mathbf A} %
\newcommand{\bfB}{\mathbf B} %
\newcommand{\bfC}{\mathbf C} %
\newcommand{\bfG}{\mathbf G} %
\newcommand{\bfH}{\mathbf H} %
\newcommand{\bfL}{\mathbf L} %
\newcommand{\bfP}{\mathbf P} %
\newcommand{\bfQ}{\mathbf Q} %
\newcommand{\bfR}{\mathbf R} %
\newcommand{\bfT}{\mathbf T} %
\newcommand{\bfS}{\mathbf S} %
\newcommand{\bfU}{\mathbf U} %
\newcommand{\bfV}{\mathbf V} %
\newcommand{\bfX}{\mathbf X} %

\DeclareMathOperator{\GL}{GL} %
\DeclareMathOperator{\U}{U} %

\newcommand\inverse{{^{-1}}}

\newcommand\GG{{\mathbb G}}
\newcommand\NN{{\mathbb N}}
\newcommand\ZZ{{\mathbb Z}}
\newcommand\FF{\mathbb F}

\newcommand\QQ{\mathbb Q}

\newcommand\CC{\mathcal C}
\newcommand\CR{\mathcal R}

\newcommand\uni{{\mathrm{uni}}}
\newcommand\uu{{\mathrm{u}}}

\usepackage{amssymb}
\usepackage{fullpage}

\numberwithin{equation}{section}

\theoremstyle{plain}
\newtheorem{lem}[equation]{Lemma}
\newtheorem{thm}[equation]{Theorem}
\newtheorem{cor}[equation]{Corollary}
\newtheorem{prop}[equation]{Proposition}

\theoremstyle{definition}
\newtheorem{exmp}[equation]{Example}

\theoremstyle{remark}
\newtheorem{rem}[equation]{Remark}

\title{On conjugacy of unipotent elements in finite groups of Lie type}

\author[S.~M.~Goodwin]
{Simon M.~Goodwin}
\address{School of Mathematics, University of Birmingham,
Birmingham, B15 2TT, United Kingdom}
\email{goodwin@maths.bham.ac.uk}
\urladdr{http://web.mat.bham.ac.uk/S.M.Goodwin}

\author[G.\ R\"ohrle]{Gerhard R\"ohrle}
\address%[G.~R\"{o}hrle]
{Fakult\"at f\"ur Mathematik, Ruhr-Universit\"at Bochum, D-44780
Bochum, Germany} \email{gerhard.roehrle@rub.de}
\urladdr{http://www.ruhr-uni-bochum.de/ffm/Lehrstuehle/Lehrstuhl-VI/rubroehrle.html}
 

\thanks{2000 {\it Mathematics Subject Classification}.
Primary 20G40, 20E45, Secondary 20D20.}

\begin{document}

\begin{abstract}
Let $\bfG$ be a connected reductive algebraic group defined over
$\FF_q$, where $q$ is a power of a prime $p$ that is good for
$\bfG$. Let $F$ be the Frobenius morphism associated with the
$\FF_q$-structure on $\bfG$ and set $G = \bfG^F$, the fixed point
subgroup of $F$.  Let $\bfP$ be an $F$-stable parabolic subgroup of
$\bfG$ and let $\bfU$ be the unipotent radical of $\bfP$; set $P =
\bfP^F$ and $U = \bfU^F$. Let $G_\uni$ be the set of unipotent
elements in $G$. In this note we show that the number of conjugacy
classes of $U$ in $G_\uni$ is given by a polynomial in $q$ with
integer coefficients.
\end{abstract}

\maketitle

\section{Introduction}
Let $\U_n(q)$ be the subgroup of $\GL_n(q)$ consisting of upper
unitriangular matrices, where $q$ is a power of a prime. A
longstanding conjecture attributed to G.~Higman (cf.\ \cite{higman})
states that the number of conjugacy classes of $\U_n(q)$ for fixed
$n$ is a polynomial in $q$ with integer coefficients. This
conjecture has been verified for $n \le 13$ by computer calculation
in work of A.~Vera-Lopez and J.~M.~Arregi, see
\cite{veralopezarregi}. There has been much interest in this
conjecture, for example from G.R.~Robinson \cite{robinson} and
J.~Thompson \cite{thompson}.

In \cite{alperin} J.~Alperin showed that a related question is
easily answered, namely that the number of $\U_n(q)$-conjugacy
classes in all of $\GL_n(q)$ is a polynomial in $q$ with integer
coefficients.   In \cite{goodwinroehrle:unipotent} the authors
generalized Alperin's result twofold, by replacing $\GL_n(q)$ by a
finite group of Lie type $G$ and by replacing $\U_n(q)$ by the
unipotent radical $U$ an arbitrary parabolic subgroup $P$ of $G$.
Precisely, in \cite[Thm.\ 4.5]{goodwinroehrle:unipotent}, under the
assumptions that the reductive algebraic group $\bfG$ corresponding
to $G$ has connected centre and that $q$ is a power of a good prime
for $\bfG$, we showed that, the number $k(U,G)$ of $U$-conjugacy
classes in $G$ is a polynomial in $q$ with integer coefficients (if
$G$ has a simple component of type $E_8$, then there exist
polynomials $m^i(z) \in \ZZ[z]$ for $i = \pm 1$ so that $k(U,G)=
m^i(q)$, when $q$ is congruent $i$ modulo $3$).

Using the machinery developed in \cite{goodwinroehrle:unipotent}, we
discuss the following related conjugacy problem in this note: we
show that the number $k(U, G_\uni)$ of $U$-conjugacy classes in the
set $G_\uni$ of unipotent elements of $G$ is a polynomial in $q$
with integer coefficients (again, if $G$ has a simple component of
type $E_8$, then two polynomials are required depending on the
congruence class of $q$ modulo $3$); see Theorem \ref{thm:main} for
a precise statement.  For this theorem we do not require the
assumption that the centre of $\bfG$ is connected; this is because
$U$ and $G_\uni$ are ``independent'' up to isomorphism of the
isogeny class of $\bfG$.

One can view Alperin's result in \cite{alperin} and Theorem
\ref{thm:main} for $G = \GL_n(q)$ and $U = \U_n(q)$ as evidence in
support of Higman's conjecture. In \cite{alperin} Alperin remarks
that it is unlikely to be possible to obtain a proof of Higman's
conjecture by descent from his theorem; it seems equally improbable
that a proof of this conjecture can be deduced from Theorem
\ref{thm:main}.

\smallskip

As general references on finite groups of Lie type, we refer the
reader to the books by Carter \cite{carter} and Digne--Michel
\cite{dignemichel}.

\section{Preliminaries}
\label{s:prelims}

\subsection{General notation for algebraic groups}
\label{s:notn} We introduce some notation used throughout.  Let $q$
be a power of a prime $p$.  By $\FF_q$ we denote the field of $q$
elements and by $\overline{\FF}_q$ its algebraic closure. Throughout
this paper, we identify algebraic groups defined over $\FF_q$ with
their group of $\overline{\FF}_q$-rational points.  So in particular, the
additive group $\GG_a$ and the multiplicative group $\GG_m$ are
identified with the additive group $\overline{\FF}_q$ and
multiplicative group $\overline{\FF}_q^*$ respectively.

Let $\bfG$ be a connected reductive algebraic group defined over
$\FF_q$, where $p$ is assumed to be good for $\bfG$. Let $F$ be the
Frobenius morphism associated with the $\FF_q$-structure on $\bfG$
and set $G = \bfG^F$ the finite group of fixed points of $F$ in
$\bfG$.

Let $\bfH$ be a closed $F$-stable subgroup of $\bfG$.  We write
$\bfH^\circ$ for the identity component of $\bfH$, $\bfH_\uni$ for
the subset of unipotent elements in $\bfH$ and $H = \bfH^F$.  By
$|H|_p$ we denote the size of a Sylow $p$-subgroup of $H$ and by
$|H|_{p'}$ the $p'$-part of the order of $H$.  Let $S$ be an
$H$-stable subset of $G$.  We write $k(H,S)$ for the number of
$H$-conjugacy classes in $S$.  Given $x \in G$ we write
$C_{\bfH}(x)$ for the centralizer of $x$ in $\bfH$ and $C_H(x)$ for
the centralizer of $x$ in $H$; for $x \in H$, then we write $C_S(x)$
for the set of fixed points of $x$ in $S$. The $H$-conjugacy class
of $x$ is denoted by $H \cdot x$. We write
\begin{equation*}
f_H^G(x) = |\{{}^gH \mid x \in {}^gH, g \in G\}|
\end{equation*}
for the number of conjugates of $H$ in $G$ containing $x$.

\subsection{Axiomatic setup for connected reductive algebraic groups}
\label{s:setup}

For the statement of our main theorem (Theorem \ref{thm:main}) we
require the axiomatic setup for connected reductive algebraic groups
given in \cite[\S 2.2]{goodwinroehrle:unipotent}, which we now
recall for completeness and convenience. The idea is that a tuple of
combinatorial objects is used to define a family of connected
reductive groups indexed by prime powers. We refer the reader to
\cite[\S 0, \S 3]{dignemichel} for some of the results used below.

Let $\Psi = (X, \Phi, \check X, \check\Phi)$ be a root datum.  Then
given a finite field $\FF_q$, the root datum $\Psi$ determines a
connected reductive algebraic group $\bfG$ over $\overline{\FF}_q$
and a maximal torus $\bfT$ of $\bfG$ such that $\Psi$ is
the root datum of $\bfG$ with respect to $\bfT$.
Let $\Pi$ be a base for $\Phi$; this determines a Borel subgroup
$\bfB$ of $\bfG$ containing $\bfT$.

Let $F_0 : X \to X$ be an automorphism of finite order such that
$F_0(\Phi) = \Phi$, $F_0(\Pi) = \Pi$ and $F_0^*(\check\Phi) =
\check\Phi$.  Then for any prime power $q$, the automorphism $F_0$
defines a Frobenius morphism $F: \bfG \to \bfG$ such that the induced
action of $F$ on $X$ is given by $q \cdot F_0$.  Further, $\bfB$ and
$\bfT$ are $F$-stable, so that $\bfT$ is a maximally split maximal torus
of $\bfG$.

A subset $J$ of $\Pi$ determines the standard parabolic subgroup
$\bfP = \bfP_J$ of $\bfG$.  If $F_0(J) = J$, $q$ is a prime power
and $F$ is the corresponding Frobenius morphism of $\bfG$, then
$\bfP$ is $F$-stable.

Summing up, the discussion above implies that the quadruple $\Delta
= (\Psi,\Pi,F_0,J)$, along with a prime power $q$ determines:
\begin{itemize}
\item a connected reductive algebraic group $\bfG$ defined over $\FF_q$ with
corresponding Frobenius morphism $F$;
\item a maximally split $F$-stable maximal torus $\bfT$;
\item an $F$-stable Borel subgroup $\bfB \supseteq \bfT$ of $\bfG$; and
\item an $F$-stable parabolic subgroup $\bfP \supseteq \bfB$.
\end{itemize}

The notation we use for $\bfG$, $\bfB$, $\bfT$ and $\bfP$ does not
reflect the fact that their $\FF_q$-structure depends on the choice
of a prime power $q$. Let $q$ be a prime power and $m$ a positive
integer, write $F$ for the Frobenius morphism corresponding to $q$.
Then it is not necessarily the case that the Frobenius morphism
corresponding to the prime power $q^m$ is $F^m$, i.e.\ the
definition of $\bfG$ over $\FF_{q^m}$ is not necessarily obtained
from the $\FF_q$-structure by extending scalars.  The definitions of
$\bfG$ over $\FF_{q^m}$ are not equivalent if $F_0$ is not the
identity and there is a common divisor of $m$ and the order $F_0$.
However, in order to keep the notation short, we choose not to show
this dependence on $q$. We refer the reader to \cite[Rem.\
2.1]{goodwinroehrle:unipotent} for further explanation of our
convention for varying $q$.

Given the data $\Delta = (\Psi,\Pi,F_0,J)$ and prime power $q$, we
note that the unipotent radical $\bfU = R_\uu(\bfP)$ of $\bfP =
\bfP_J$ and the unique Levi subgroup $\bfL = \bfL_J$ of $\bfP$
containing $\bfT$ are determined. Since both $\bfP$ and $\bfT$ are
$F$-stable, so are $\bfU$ and $\bfL$.

\subsection{Commuting varieties}
\label{s:commuting}

Let $\bfH$ and $\bfS$ be a closed subgroup and a closed
$\bfH$-stable subvariety of $\bfG$, respectively. The {\em commuting
variety of $\bfH$ and $\bfS$} is the closed subvariety of $\bfH
\times \bfS$ defined by
\[
\CC(\bfH,\bfS) = \{(h,s) \in \bfH \times \bfS \mid hs = sh \}.
\]
Assume that both $\bfH$ and $\bfS$ are $F$-stable. Then $F$ acts on
$\CC(\bfH,\bfS)$ and we have $\CC(\bfH,\bfS)^F = \CC(\bfH^F,\bfS^F)
= \CC(H,S)$. The Burnside counting formula gives
\begin{equation}
\label{eq:comm} |\CC(H,S)| = \sum_{x\in H}|C_S(x)| = |H|\cdot
k(H,S).
\end{equation}

\subsection{Kempf--Rousseau theory}
\label{sub:git}

We now briefly recall the theory of optimal cocharacters from
geometric invariant theory.  We require this in the proof of
Proposition \ref{prop:centralizer}, which is key to the proof of
Theorem \ref{thm:main}.

Let $\check X(\bfG)$ denote the set of cocharacters of $\bfG$, i.e.,
the set of homomorphisms $\mathbb G_m \to \bfG$. There is a left
action of $\bfG$ on $\check X(\bfG)$: for $\mu\in \check X(\bfG)$
and $g\in \bfG$ we define $g\cdot \mu \in \check X(\bfG)$ by
$(g\cdot \mu)(t) = g\mu(t)g\inverse$.

Let $\bfX$ be an affine variety.  Let $\phi : \GG_m \to \bfX$ be a
morphism of algebraic varieties. We say that $\underset{t\to
0}{\lim}\, \phi(t)$ exists if there exists a morphism $\widehat\phi
:\GG_a \to \bfX$ (necessarily unique) whose restriction to $\GG_m$
is $\phi$; if this limit exists, then we set $\underset{t\to
0}{\lim}\, \phi(t) = \widehat\phi(0)$.

Let $\bfP$ be a parabolic subgroup of $\bfG$, and let $\bfL$ be a
Levi subgroup of $\bfP$.  We recall, see for example \cite[Prop.\
8.4.5]{spr2}, that there exists $\lambda \in \check X(\bfG)$ such
that: $\bfP = \bfP_\lambda := \{g\in \bfG \mid \underset{t\to
0}{\lim}\, \lambda(t) g \lambda(t)^{-1} \textrm{ exists}\}$; $\bfL =
\bfL_\lambda := C_\bfG(\lambda(\GG_m))$; and $R_\uu(\bfP) = \{g\in
\bfG \mid \underset{t\to 0}{\lim}\, \lambda(t) g \lambda(t)^{-1} =
1\}$. Moreover, the map $c_\lambda : \bfP_\lambda \to \bfL_\lambda$
given by $g \mapsto \underset{t\to 0}{\lim}\,
\lambda(t)g\lambda(t)^{-1}$ is a homomorphism of algebraic groups
with kernel $\ker c_\lambda = R_\uu(\bfP)$; we note that $c_\lambda$
is simply the projection from $\bfP$ onto $\bfL$ along the
semidirect decomposition $\bfP = \bfL R_\uu(\bfP)$.

Let $\bfG$ act on the affine variety $\bfX$. For $x \in \bfX$ let
$\bfG \cdot x$ denote the $\bfG$-orbit of $x$ in $\bfX$ and
$C_\bfG(x)$ the stabilizer of $x$ in $\bfG$. Let $x \in \bfX$ and
let $\bfC$ be the unique closed orbit in the closure of $\bfG\cdot
x$, we refer the reader to \cite[1.3]{rich} for a proof that there
is a unique closed $\bfG$-orbit in $\overline{\bfG \cdot x}$. The
Kempf--Rousseau theory tells us that there exists a non-empty subset
$\Omega(x)$ of $\check X(\bfG)$ consisting of so called
\emph{optimal cocharacters} $\lambda$ such that $\underset{t\to
0}{\lim}\,\lambda(t) \cdot x$ exists and belongs to $\bfC$, we refer
the reader to \cite{kempf} or \cite{rousseau} for information on the
Kempf--Rousseau theory and the definition of optimal cocharacters.
Moreover, there exists a parabolic subgroup $\bfP(x)$ of $\bfG$ so
that $\bfP(x) = \bfP_\lambda$ for every $\lambda \in \Omega(x)$, and
we have that $\Omega(x)$ is a single $\bfP(x)$-orbit. Further, for
every $g \in \bfG$, we have $\Omega(g\cdot x) = g\cdot \Omega(x)$
and $\bfP(g \cdot x) = g\bfP(x)g\inverse$. In particular, $C_\bfG(x)
\le \bfP(x)$. The parabolic subgroup $\bfP(x)$ is called the
\emph{optimal} or \emph{destabilizing parabolic subgroup associated
to $x$}.

\begin{rem}
\label{rem:F-Kempf} Suppose that $\bfG$, $\bfX$ and the action of
$\bfG$ on $\bfX$ are all defined over $\FF_q$ and let $F$ denote the
Frobenius morphism associated with the $\FF_q$-structures on both
$\bfG$ and $\bfX$. There is an action of $F$ on $\check X(\bfG)$ as
follows: for $\mu\in \check X(\bfG)$ we define $F \cdot \mu \in
\check X(\bfG)$ by $(F\cdot \mu)(t) = F(\mu(F^{-1}(t)))$, where $F :
\GG_m \to \GG_m$ is given by  $F(t) = t^q$, see \cite[\S 4]{kempf}.
Thanks to \cite[Thm.~4.2]{kempf} and \cite[\S 2]{LMS}, if $x$ is
fixed by $F$, then both $\Omega(x)$ and $\bfP(x)$ are $F$-stable.
\end{rem}

\section{Polynomial behaviour of $k(U,G_\uni)$}

We maintain the notation and assumptions made in the previous
sections. In particular, $\bfG$ is a connected reductive algebraic
group defined over $\FF_q$, where $q$ is a power of the prime $p$
which is good for $\bfG$.

We begin by stating a counting lemma for finite groups from
\cite{alperin}, see also \cite[Lem.\ 4.1]{goodwinroehrle:unipotent};
the argument used to prove \cite[Lem.\
4.1]{goodwinroehrle:unipotent}, which uses the Burnside counting
lemma, easily generalizes to the present situation, so we do not
include it here.

\begin{lem}
\label{lem1}
Let $\bfP = \bfL \bfU$ be an $F$-stable parabolic subgroup of $\bfG$.
Then the number of $U$-conjugacy classes in $G_\uni$ is given by
\begin{equation*}
k(U,G_\uni) = |L|\sum_{x \in \CR} \frac{|C_G(x)_\uni|}{|C_G(x)|} f_U^G(x),
\end{equation*}
where $\CR$ is a set of representatives of the unipotent $G$-conjugacy classes.
\end{lem}

Armed with the theory of optimal cocharacters from \S \ref{sub:git},
we are able to provide the following key result for our proof that
$k(U,G_\uni)$ is a polynomial in $q$.  We note that the Levi
decomposition of $C_\bfG(u)$ stated in Proposition
\ref{prop:centralizer}(i) is well-known, see for example \cite[Thm.\
A]{premet}.

\begin{prop}
\label{prop:centralizer}
Let $u \in G_\uni$.
Then
\begin{itemize}
\item[(i)]
$C_\bfG(u)$ admits a Levi decomposition $C_\bfG(u) = \bfC(u)
\bfR(u)$ with $\bfC(u) \cap \bfR(u) = \{1\}$, $\bfC(u)$ reductive
and $\bfR(u)$ the unipotent radical of $C_\bfG(u)$, such that
$\bfC(u)$ is $F$-stable; therefore, setting $C(u) = \bfC(u)^F$ and
$R(u) = \bfR(u)^F$, we obtain a Levi decomposition $C_G(u) =
C(u)R(u)$ of $C_G(u)$;
\item[(ii)]
both $C_\bfG(u)_\uni$ and $C_G(u)_\uni$ admit a ``Levi
decomposition'', $C_\bfG(u)_\uni = \bfC(u)_\uni \bfR(u)$ and
$C_G(u)_\uni = C(u)_\uni R(u)$.
\end{itemize}
\end{prop}

\begin{proof}
(i). Since $p$ is good for $\bfG$, the centralizer of $u$ in $\bfG$
has a Levi decomposition, $C_\bfG(u) = \bfC(u) \bfR(u)$ with
$\bfC(u)$ reductive and $\bfR(u)$ the unipotent radical of
$C_\bfG(u)$, see for example \cite[Thm.\ A]{premet}.  More
precisely, by \cite[Thm.\ 2.1; Prop.\ 2.5]{premet}, let $\bfP(u)$ be
the destabilizing parabolic subgroup associated to $u$, then as
explained in \S \ref{sub:git}, we have $C_\bfG(u) \subseteq
\bfP(u)$. Let $\lambda \in \Omega(u)$ be an optimal cocharacter of
$\bfG$ associated to $u$. We have $\bfP(u) = \bfP_\lambda =
\bfL_\lambda\bfU(u)$, where $\bfU(u) = R_\uu(\bfP(u))$. Setting
$\bfC(u) = \bfL_\lambda \cap C_\bfG(u)$ and $\bfR(u) = \bfU(u) \cap
C_\bfG(u)$, we obtain a Levi decomposition $C_\bfG(u) = \bfC(u)
\bfR(u)$ of $C_\bfG(u)$.

By Remark \ref{rem:F-Kempf}, both $\Omega(u)$ and $\bfP(u)$ are
$F$-stable, since $u$ $\in G$. Further, $\Omega(u)$ is a single
$\bfP(u)$-orbit. Since $\bfP(u)$ is connected and $F$-stable, it
follows, from for example \cite[Cor.\ 3.12]{dignemichel}, that there
exists an $F$-stable cocharacter in $\Omega(u)$.  We may therefore
assume that $\lambda$ is $F$-stable.  Then $\bfC(u)$ is $F$-stable.
Clearly, $\bfR(u)$ is also $F$-stable. Since $\bfC(u) \cap \bfR(u) =
\{1\}$, it follows that $C_G(u) = C_\bfG(u)^F = \bfC(u)^F \bfR(u)^F
= C(u) R(u)$.

(ii). Let $v \in C_\bfG(u)$ be unipotent. Thanks to the Levi
decomposition $C_\bfG(u) = \bfC(u) \bfR(u) \subseteq \bfP(u)$ of
$C_\bfG(u)$ from part (i), we have $v = xy$ with $x \in \bfC(u)$ and
$y \in \bfR(u)$.  We have $x = c_\lambda(v)$ where $c_\lambda :
\bfP(u) \to \bfL_\lambda$ is the canonical homomorphism defined in
\S \ref{sub:git}.  Therefore, $x$ is unipotent and we obtain the
decomposition $C_\bfG(u)_\uni = \bfC(u)_\uni \bfR(u)$.

Since $\bfC(u)_\uni \cap \bfR(u) = \{1\}$, it follows that
$C_G(u)_\uni = C_\bfG(u)_\uni^F = \bfC(u)_\uni^F \bfR(u)^F =
C(u)_\uni R(u)$, as desired.
\end{proof}

Next we require a result regarding the independence in $q$ of the
orders of centralizers of unipotent elements in $G$; more precisely
that these orders are given by a polynomial in $q$. We use the
axiomatic setup from \S \ref{s:setup} to achieve this. Fix
$(\Psi,\Pi,F_0)$, where $\Psi = (X, \Phi, \check X, \check\Phi)$,
and for a prime power $q$, let $\bfG$ and $F$  be the connected
reductive group and Frobenius morphism determined by
$(\Psi,\Pi,F_0)$ and $q$.  We assume that $X/\ZZ\Phi$ is torsion
free, where $\ZZ\Phi$ denotes the root lattice of $\bfG$, this
ensures that the centre of $G$ is connected for all $q$.  We also
assume that $q$ is a power of a good prime for $\bfG$.

Under these assumptions the parametrization of the unipotent
conjugacy classes of $G$ is independent of $q$, see for example
\cite[Prop.\ 2.5]{goodwinroehrle:unipotent}.  We let $\CR$ be a set
of representatives of the unipotent conjugacy classes of $G$, and we
use the convention of \cite[Rem.\ 2.6]{goodwinroehrle:unipotent} to
vary $u \in \CR$ with $q$.  With these conventions we can state and
prove the following proposition, which is crucial for our proof of
Theorem \ref{thm:main}.

\begin{prop} \label{P:centsize}
Assume that $X/\ZZ\Phi$ is torsion free and that $q$ is a power
of a good prime for $\bfG$. Let $u \in \CR$.  Then the order of
$C_G(u)$ is given by a polynomial in $q$.  Further, the order of
$C_G(u)_\uni$ is given by a fixed power of $q$.
\end{prop}

\begin{proof}
That the order of $C_G(u)$ is a polynomial in $q$, can been seen
from the Lusztig--Shoji algorithm for computing Green functions, see
\cite{lusztig} and \cite{shoji}. It is straightforward to see that
the order of the centralizers of unipotent elements of $G$ can be
determined from the block-diagonal matrix $\Lambda$, defined in
\cite[\S 5]{shoji}; the blocks are determined by the Springer
correspondence.
The (unknown) matrix $\Lambda$ satisfies the equation
\begin{equation} \label{e:Lambda}
{}^t P \Lambda P = \Pi,
\end{equation}
where $P$ is an (unknown) upper triangular block matrix with each
diagonal block a matrix with entries in $\QQ$, and $\Pi$ is a known
matrix with entries that are rational functions in $q$ with
coefficients independent of $q$, see \cite[(5.6)]{shoji}. As stated
in {\em loc.\ cit.}, the matrix $\Lambda$ is uniquely determined by
\eqref{e:Lambda}; moreover, one sees that the entries of $\Lambda$
are rational functions in $q$, with coefficients independent of $q$.
In particular, we can deduce that $|C_G(u)|$ is a rational function in
$q$, and then a standard argument, see for example \cite[Lem.\
2.12]{goodwinroehrle:unipotent}, tells us that $|C_G(u)|$ is in fact
a polynomial in $q$.

The second statement in the lemma now follows from Steinberg's
formula applied to the Levi factor $C(u)$ of $C_G(u)$, see for
example \cite[Cor.\ 9.5]{dignemichel}.
\end{proof}

\begin{rem}
It seems likely that a stronger result than Proposition
\ref{P:centsize} regarding the structure of $C_G(u)$ holds. That is
the root datum corresponding to $\bfC(u)^\circ$, the component group
$\bfA(u)$ of $\bfC(u)$, and the action of $F$ on the root lattice
for $\bfC(u)^\circ$ and $\bfA(u)$ do not depend on $q$. It is known
that the root datum of $\bfC(u)$ does not depend on $q$, though this
is only by a case by case analysis, see for example the discussion
at the end of \cite[5.11]{jantzen}. There is a general proof that
the structure of $\bfA(u)$ does not depend on $q$, see
\cite{mcninchsommers} or \cite{premet}.  One then needs to check
that the action of $F$ on $\bfC(u)$ for {\em split elements} $u$ is
independent of $q$; for $G$ of type $E_8$ one other case needs to be
dealt with separately. Further one needs to know that the action of
$\bfA(u)$ on the set of simple roots of $\bfC(u)^\circ$ does not
depend on $q$.  We have chosen not to pursue this here.
\end{rem}

We are now in a position to prove the principal result of this
paper, which is an analogue of \cite[Thm.\
4.5]{goodwinroehrle:unipotent}.  We continue to use the axiomatic
setup from \S \ref{s:setup}.

\begin{thm}
\label{thm:main} Fix the data $\Delta = (\Psi,\Pi,F_0,J)$, where
$\Psi = (X, \Phi, \check X, \check\Phi)$. For a prime power $q$, let
$\bfG$, $F$ and $\bfP$ be the connected reductive group, Frobenius
morphism and $F$-stable parabolic subgroup of $\bfG$ determined by
$\Delta$ and let $\bfU = R_\uu(\bfP)$. Assume that $q$ is power
of a good prime for $\bfG$.
\begin{itemize}
\item[(i)]
Suppose that $\bfG$ does not have a simple component of type $E_8$.
Then there exists $m(z) \in \ZZ[z]$ such that $k(U,G_\uni) = m(q)$.
\item[(ii)]
Suppose that $\bfG$ has a simple component of type $E_8$. Then there exist
$m^i(z) \in \ZZ[z]$ ($i = \pm 1$), such that $k(U,G_\uni) = m^i(q)$,
when $q$ is congruent to $i$ modulo $3$.
\end{itemize}
\end{thm}

\begin{proof}
We begin by assuming that $X/\ZZ\Phi$ is torsion free, so that the
centre of $\bfG$ is connected.  We write $\bfL$ for the Levi
subgroup of $\bfP$ containing $\bfT$.

By Lemma \ref{lem1}, we have
\begin{equation} \label{e:k(U,G_uni)}
k(U,G_\uni) = |L|\sum_{x \in \CR} \frac{|C_G(x)_\uni|}{|C_G(x)|} f_U^G(x),
\end{equation}
where $\CR$ is a set of representatives of the unipotent
$G$-classes. With the assumptions that $X/\ZZ\Phi$ is torsion free
and that $q$ is power of a good prime for $\bfG$, it follows from
\cite[Prop.\ 2.5]{goodwinroehrle:unipotent} that the set $\CR$ is
independent of $q$, where we use the convention of \cite[Rem.\
2.6]{goodwinroehrle:unipotent} to vary $x$ with $q$.

Since $L$ is a finite reductive group, the factor $|L|$ is a
polynomial in $q$ (\cite[p.\ 75]{carter}). Thanks to \cite[Lem.\
3.1(ii), Thm.\ 3.10]{goodwinroehrle:unipotent}, each of the factors
$f_{U}^G(x)$ in the sum above is a polynomial in $q$, unless we are
in case (ii) when $f_{U}^G(x)$ is given by two polynomials depending
on $q$ modulo $3$.  From Proposition \ref{P:centsize} we have that
$|C_G(x)|$ and $|C_G(x)_\uni|$ are polynomials in $q$.  Hence,
$k(U,G_\uni)$ is a rational function in $q$.  Now by a standard
argument, see for example \cite[Lem.\
2.12]{goodwinroehrle:unipotent}, we can conclude that $k(U,G_\uni)$
is a polynomial function in $q$ with rational coefficients.

Assume now that $\bfG$ is split over $\FF_q$, i.e.\ that $F_0$ is
the identity.  Thanks to \eqref{eq:comm}, $|\CC(\bfU,\bfG_\uni)^F|$
is a polynomial in $q$ with rational coefficients.  The assumption
that $\bfG$ is split means that $|\CC(U,G_\uni)|$ gives the number
of $\FF_q$-rational points in the variety $\CC(\bfU,\bfG_\uni)$ viewed as
a variety defined over $\FF_p$.  Now using the Grothendieck trace
formula (see \cite[Thm.\ 10.4]{dignemichel}), one can prove that the
coefficients of this polynomial are integers, see for example
\cite[Prop.\ 6.1]{Re}.

A further standard argument using the Grothendieck trace formula now
tells us that the eigenvalues of $F$ on the $l$-adic cohomology groups of
$\CC(\bfU,\bfG_\uni)$ are all powers of $q$, see for example the proof of
\cite[Prop.\ 6.1]{Re}.  Now assume that $\bfG$ is not split and let
$d$ be the order of $F_0$.  Then arguments like those used to prove
\cite[Prop.\ 3.20]{goodwinroehrle:unipotent} imply that the
eigenvalues of $F$ on the $l$-adic cohomology groups of
$\CC(\bfU,\bfG_\uni)$ (viewed as a variety defined over $\FF_q$) are of the
form $\zeta q^m$, where $\zeta$ is a $d$th root of unity.  Following
the arguments to prove \cite[Prop.\ 3.20]{goodwinroehrle:unipotent},
one can now show that the coefficients of the polynomial
$|\CC(U,G_\uni)|$ are integers.
Then using \eqref{eq:comm} again, it follows that $k(U,G_\uni)$ is a
polynomial function in $q$ with integer coefficients.

\smallskip

Now remove the assumption that $X/\ZZ\Phi$ is torsion free.  Let
$\sigma : \bfG \to \hat{\bfG}$ be an isogeny that is defined over
$\FF_q$, where $\hat \bfG$ is a reductive group defined over $\FF_q$
with connected centre. Then $\sigma$ induces an isomorphism between
$\bfU$ and $\hat{\bfU}$ and between $\bfG_\uni$ and $\hat
\bfG_\uni$, since $Z(\bfG) \cap \bfU = \{1\} = Z(\bfG) \cap
\bfG_\uni$, where $Z(\bfG)$ is the centre of $\bfG$. It follows
easily that $k(U,G_\uni) = k(\hat U,\hat G_\uni)$ is given by a
polynomial in $q$ with integer coefficients.
\end{proof}

Recall that two parabolic subgroups of $\bfG$ are called
\emph{associated} if they have Levi subgroups that are conjugate in
$\bfG$. It was already remarked in \cite[Cor.\
3.5]{goodwinroehrle:unipotent} that if $\bfP$ and $\bfQ$ are
associated parabolic subgroups of $\bfG$ with unipotent radicals
$\bfU$ and $\bfV$ respectively, then the functions $f_U^G$ and
$f_V^G$ are equal on unipotent elements of $G$. Therefore, from
\eqref{e:k(U,G_uni)} we can observe the following corollary in the
same way as \cite[Cor.\ 4.7]{goodwinroehrle:unipotent}.

\begin{cor}
\label{cor2} Let $\bfP$ and $\bfQ$ be associated $F$-stable
parabolic subgroups of $\bfG$ with unipotent radicals $\bfU$ and
$\bfV$ respectively. Then
\[
k(U,G_\uni) = k(V,G_\uni).
\]
\end{cor}

In order to be able to compute the polynomials given in Theorem
\ref{thm:main} explicitly, we reformulate the expression for
$k(U,G_\uni)$ in \eqref{e:k(U,G_uni)} in terms of Green functions.
Using \cite[Lem.\ 3.3]{goodwinroehrle:unipotent}, Proposition
\ref{prop:centralizer}, and Steinberg's formula \cite[Cor.\
9.5]{dignemichel}, we obtain
\begin{align}
\label{eqn:eq2}
k(U,G_\uni) & \notag = |L|\sum_{x \in \CR} \frac{|C_G(x)_\uni|}{|C_G(x)|}
 \left(\frac{1}{|L|_{p}|W_\bfL|}\sum_{w\in W_\bfL}(-1)^{l(w)}Q_{\bfT_w}^\bfG(x) \right)\\
 & = \frac{1}{|W_\bfL|}\frac{|L|\ }{|L|_{p}}\;\sum_{x \in \CR}\frac{|C(x)|_p^2}{|C(x)|}
           \left(\sum_{w \in W_\bfL} (-1)^{l(w)} Q_{\bfT_w}^\bfG(x)\right)\\
 & = \notag \frac{1}{|W_\bfL|}|L|_{p'}\sum_{x \in \CR}\frac{|C(x)|_p}{|C(x)|_{p'}}
           \left(\sum_{w \in W_\bfL} (-1)^{l(w)} Q_{\bfT_w}^\bfG(x)\right),
\end{align}
where $C(x)$ is the reductive part of $C_G(x)$, as in Proposition
\ref{prop:centralizer}(i), $W_\bfL$ is the Weyl group of $\bfL$; for
$w \in W_\bfL$ we write $\bfT_w$ for the $F$-stable twisted torus
associated with $w \in W_\bfL$, and $Q_{\bfT_w}^\bfG$ is the Green
function associated with $\bfT_w$. For more information on Green
functions we refer the reader to \cite[\S 7.6]{carter}.
We note that the sum in \eqref{eqn:eq2} is effectively only over
representatives of the $\bfG$-orbits that meet $\bfU$.  This follows
from the fact that the term $f_U^G(x)$ in \eqref{e:k(U,G_uni)} is
obviously zero if $\bfG \cdot x \cap \bfU = \varnothing$.

Using the {\tt chevie} package in GAP3 (\cite{gap}) along with some
code provided by M.~Geck, it is possible to explicitly calculate the
polynomials $m(z)$ in Theorem \ref{thm:main} for $\bfG$ of small
rank. We illustrate this with some examples for the case $\bfG =
\GL_n$ and $\bfP = \bfB$ is a Borel subgroup of $\bfG$.

\begin{exmp}
\label{ex:polys2} In Table \ref{Tab:GL} below we give the
polynomials for $k(U,G_\uni)$ in case $G = \GL_n(q)$ and $P = B$ is
a Borel subgroup of $G$, for $n = 2,\dots,10$.  In this case we take
$U = \U_n(q)$ to be the group of upper unitriangular matrices.

\begin{table}[h!tb]
\renewcommand{\arraystretch}{1.5}
\begin{tabular}{|l|p{400pt}|}
\hline $n$ & $k(\U_n(q),\GL_n(q)_\uni)$ \\
\hline\hline 2 &  $2q - 1$ \\
\hline
3 & $q^3 + 3q^2 - 3q$ \\
\hline
4 & $q^6 + 5q^4 - 9q^2 + 4q$  \\
\hline 5 & $q^{10} + 4q^7 + 4q^6 + 6q^5 - 20q^4 - 10q^3 + 21q^2 - 4q
-
1$ \\
\hline 6 & $q^{15} + 5q^{11} - q^{10} + 13q^9 - 5q^8 + q^7 - 20q^6 -
44q^5 + 52q^4 + 25q^ 3 - 31q^2 + 5q$
 \\
\hline 7 & $q^{21} + 6q^{16} - q^{15} - q^{14} + 13q^{13} + 22q^{12}
- 41q^{11} + 37q^{10} - 49q^ 9 - 57q^8 - 12q^7 + 71q^6 + 139q^5 -
120q^4 - 51q^3 +
50q^2 - 5q - 1$ \\
\hline 8 & $q^{28} + 7q^{22} - q^{21} - q^{20} - q^{19} + 19q^{18} +
13q^{17} + 7q^{16} - 54q^{15} + 39q^{14} + 39q^{13} - 125q^{12} -
86q^{11} - 93q^{10} + 225q^9 + 160q^8 + 102q^7 - 164q^6 - 322q^5 +
207q^4 + 87q^3 - 64q^2
+ 6q$  \\
\hline 9 & $q^{36} + 8q^{29} - q^{28} - q^{27} - q^{26} - q^{25} +
26q^{24} + 19q^{23} - 44q^{22} + 41q^{21} + 9q^{20} + 25q^{19} -
119q^{18} + 57q^{17} - 134q^{16} + 119q^{15} - 458q^{14} + 177q^{13}
+ 290q^{12} - 121q^{11} + 1315q^{10} - 807q^9 - 971q^8 + 60q^
7 + 326q^6 + 568q^5 - 319q^4 - 145q^3 + 89q^2 - 7q$ \\
\hline 10 &  $q^{45} + 9q^{37} - q^{36} - q^{35} - q^{34} - q^{33} -
q^{32} + 34q^{31} + 26q^{30} - 54q^{29} - 8q^{28} + 68q^{27} -
41q^{26} + 251q^{25} - 258q^{24} - 395q^{23} + 474q^{22} + 259q^{21}
- 674q^{20} - 174q^{19} + 126q^{18} - 1384q^{17} + 3300q^{16} -
1299q^{15} - 1227q^{14} + 4050q^{13} - 2400q^{12} - 691q^{11} -
2676q^{10} + 944q^9 + 3298q^8 - 808q^7 - 293q^6 - 1017q^5 + 455q^4 +
210q^3 - 108q^2 + 8q$
\\
\hline
\end{tabular}
\medskip
\caption{$k(\U_n(q),\GL_n(q)_\uni)$} \label{Tab:GL}
\end{table}
\end{exmp}

\begin{rem}
\label{q-1} One can easily check that each polynomial
$k(\U_n(q),\GL_n(q)_\uni)$  in Table \ref{Tab:GL} when expressed as
a polynomial in $q-1$ has all coefficients non-negative. One might
conjecture that indeed in general each of the polynomials satisfies
$k(U,G_\uni) \in \NN[q-1]$.  As stated in \cite[Rem.\
4.13]{goodwinroehrle:unipotent}, this is also the case for each of
the explicit examples of the polynomials $k(U,G)$ calculated in
\cite{goodwinroehrle:unipotent}. It would be interesting to know if
there is a geometric explanation for these positivity phenomena.
\end{rem}

\begin{rem}
Let $\bfP$ be an $F$-stable parabolic subgroup of $\bfG$. Clearly,
we have
\begin{equation*}
\label{eq:k1}
k(U, G_\uni) = \sum_{u \in \CR}k(U, G\cdot u),
\end{equation*}
where $\CR$ is a complete set of representatives of the unipotent
$G$-conjugacy classes. By an analogue of Lemma \ref{lem1}, we get
\begin{equation*}
\label{eq:k2}
k(U,G_\uni) = |L|\sum_{u \in \CR} \left(\sum_{x \in \CR}
\frac{|C_G(x) \cap G\cdot u|}{|C_G(x)|} f_U^G(x)\right).
\end{equation*}
It would be interesting to know whether each of the summands $k(U,
G\cdot u)$ is a polynomial in $q$; this is the case if $|C_G(x) \cap
G\cdot u|$ is a polynomial in $q$ for all $x$ and $u$.
\end{rem}

\begin{rem}
Using arguments as in \cite{goodwinroehrle:parabolic}, it is
possible to show that in case $\bfG = \GL_n$, the number of
$P$-conjugacy classes in $G_\uni$ is given by a polynomial in $q$.
As the details are technical, we choose not to include them here.
For arbitrary $G$ and $P$ it is not clear whether $k(P,G_\uni)$ is
polynomial or even given by Polynomials On Residue Classes (PORC).
\end{rem}

\smallskip

{\bf Acknowledgments}: This research was funded in part by EPSRC
grant EP/D502381/1.  We are very grateful to M.~Geck for some useful
conversations about this work, and for providing the computer code used
in our examples.

%%%%%%%%%%%%%%%%%%%%%%%%%%%%%%%%%%%%%%%%%%%%%%%%%%%%%%%%%%%%%%%%%%%%%%
%%%%%%%%%%%%% bibliography
%%%%%%%%%%%%%%%%%%%%%%%%%%%%%%%%%%%%%%%%%%%%%%%%%%%%%%%%%%%%%%%%%%%%%%
%\smallskip


\begin{thebibliography}{88}

\bibitem{alperin}
J.~L.~Alperin,
\emph{Unipotent conjugacy in general linear groups},
Comm. Algebra \textbf{34} (2006), no.\ 3, 889--891.

\bibitem{carter}
R.~W.~Carter, \emph{Finite groups of Lie type. Conjugacy classes and
complex characters}, Pure and Applied Mathematics, New York, 1985.

\bibitem{dignemichel}
F.~Digne and J.~Michel, \emph{Representations of finite groups of
Lie type}, London Mathematical Society Student Texts \textbf{21},
Cambridge University Press, Cambridge, 1991.

\bibitem{gap}
The GAP group, \emph{GAP -- Groups, Algorithms, and Programming --
version 3 release 4 patchlevel 4}, Lehrstuhl D f\"ur Mathematik,
Rheinisch Westf\"alische Technische Hochschule, Aachen, Germany,
1997.

\bibitem{goodwinroehrle:unipotent}
S.~M.~Goodwin and G.~R\"ohrle, \emph{Rational points on generalized
flag varieties and unipotent conjugacy in finite groups of Lie
type}, to appear in Trans.\ Amer.\ Math.\ Soc.\ (2008).

\bibitem{goodwinroehrle:parabolic}
\bysame, %S.~M.~Goodwin and G.~R\"ohrle,
\emph{Parabolic conjugacy in general linear groups},
Journal of Algebraic Combinatorics, \textbf{27}, (2008), no. 1, 99--111.

\bibitem{higman}
G.~Higman,
\emph{Enumerating $p$-groups. I. Inequalities}, Proc.\
London Math.\ Soc.\ (3) \textbf{10} (1960) 24--30.

\bibitem{jantzen}
J.~C.~Jantzen,
emph{Nilpotent orbits in representation theory},
Lie Theory. Lie Algebras and Representations. Progress in Math.\
vol.\ 228, Birkh\"auser Boston, 2004.

\bibitem{kempf}
G.~R.~Kempf, \emph{Instability in invariant theory}, Ann.\  Math.\
\textbf{108} (1978), 299--316.

\bibitem{LMS}
M.~W.~Liebeck, B.~M.~S.~Martin and A.~Shalev, \emph{On conjugacy
classes of maximal subgroups of finite simple groups, and a related
zeta function}, Duke Math. J., \textbf{128} (2005), no. 3, 541--557.

\bibitem{mcninchsommers}
G.~McNinch, E.~Sommers, \emph{Component groups of unipotent
centralizers in good characteristic}, J.\ Algebra \textbf{260}
(2003), no. 1, 323--337.

%\bibitem{lusztig}
%G.~Lusztig,
%\emph{Character sheaves}, Adv. Math. \textbf{56} (1985), 193--237; II, \textbf{57} (1985), 226--265; III,
%\textbf{57} (1985), 266--315; IV, \textbf{59} (1986), 1-63; V, \textbf{61} (1986), 103--155.

\bibitem{lusztig}
G.~Lusztig, \emph{Character sheaves II}, Adv. Math.  \textbf{57}
(1985), 226--265; III, \textbf{57} (1985.

\bibitem{premet}
A.~Premet, \emph{Nilpotent orbits in good characteristic and the
Kempf--Rousseau theory}, J. Algebra. \textbf{260}, (2003), 338--366.

\bibitem{Re}
M.~Reineke,
\emph{Counting rational points of quiver moduli}, Int.\
Math.\ Res.\ Not.\ (2006), Art.\ ID 70456.

\bibitem{rich}
R.~W.~Richardson, {\em Affine coset spaces of reductive algebraic
groups}, Bull.\ London Math.\ Soc.\ {\bf 9} (1977), no.\ 1, 38--41.

\bibitem{rousseau}
G.~Rousseau, \emph{Immeubles sph\'eriques et th\'eorie des
invariants}, C. R. Acad. Sci. Paris \textbf{286} (1978), 247--250.

\bibitem{robinson}
G.~R.~Robinson, \emph{Counting conjugacy classes of unitriangular
groups associated to finite-dimensional algebras}, J.\ Group Theory
\textbf{1} (1998), no.\ 3, 271--274.

\bibitem{shoji}
T.~Shoji, \emph{Green functions of reductive groups over a finite
field}, The Arcata Conference on Representations of Finite Groups
(Arcata, Calif., 1986), 289--301, Proc.\ Sympos.\ Pure Math., 47,
Part 1, Amer.\ Math.\ Soc., Providence, RI, 1987.

\bibitem{spr2}
T.~A.~Springer, \emph{Linear algebraic groups}, Second edition.
Progress in Mathematics, 9. Birkh\"auser Boston, Inc., Boston, MA,
1998.

\bibitem{thompson}
J.~Thompson, \emph{$k(\U_n(F_q))$}, Preprint,
{\tt{http://www.math.ufl.edu/fac/thompson.html}}.

\bibitem{veralopezarregi}
A.~Vera-L\'opez and J.~M.~Arregi, \emph{Conjugacy classes in unitriangular
matrices}, Linear Algebra Appl.\  \textbf{370}
(2003), 85--124.

\end{thebibliography}
\end{document}